\documentclass[10pt,journal,twocolumn]{IEEEtran}

\usepackage{longtable}
\usepackage{amsmath}
\usepackage{graphicx}
\usepackage{color}
\usepackage[dvips]{epsfig}

\usepackage{graphicx, amssymb}
\usepackage{amsmath}
\newtheorem{theo}{Theorem}

\newtheorem{fact}{Fact}
\newtheorem{rem}{Remark}

\newtheorem{lem}{Lemma}

\newtheorem{pro}{Property}

\def\L2{{\cal L}_2}

\newlength{\defbaselineskip}
\newcommand{\setlinespacing}[1]%
           {\setlength{\baselineskip}{#1 \defbaselineskip}}

\newcommand{\beq}{\begin{equation}}
\newcommand{\eeq}{\end{equation}}

\newcommand{\carre} {\hfill $\blacksquare$}

\begin{document}

\title{\LARGE \bf
Synthesis of  Output-Feedback Controllers for Mixed Traffic Systems in Presence of Disturbances and Uncertainties
}

\author{Shima Sadat Mousavi$^{\dagger}$, Somayeh Bahrami$^{\ddagger}$, and Anastasios Kouvelas$^{\dagger}$% <-this % stops a space
\thanks{$^{\dagger}$The authors are with the Institute for Transport Planning and Systems, ETH Zurich, Switzerland, \texttt{\small \{shimaossadat.mousavi, anastasios.kouvelas\}@ivt.baug.ethz.ch}}%
\thanks{$^{\ddagger}$The author is with the Department of Electrical Engineering, Razi University, Kermanshah, Iran, \texttt{\small s.bahrami@razi.ac.ir}}%
}

\maketitle

%-----------------------------------------------------------------------------
{

\begin{abstract}
In this paper, we study  mixed traffic systems that move along a single-lane ring-road or open-road. The traffic flow forms a platoon, which includes a number of heterogeneous human-driven vehicles (HDVs) together with only one connected and automated vehicle (CAV) that receives information from several neighbors. The dynamics of HDVs are assumed to follow the optimal velocity model (OVM), and the acceleration of the single CAV is directly controlled by a dynamical output-feedback controller. The ultimate goal of this work is to present a robust control strategy that can smoothen the traffic flow in the presence of undesired disturbances (e.g.\ abrupt deceleration) and parametric uncertainties. A prerequisite for synthesizing a dynamical output controller is the stabilizability and detectability of the underlying system.
Accordingly, a theoretical analysis is presented first to prove the stabilizability and detectability of the mixed traffic flow system. Then, two $H_\infty$  control strategies, with and without considering uncertainties in the system dynamics, are designed. The efficiency of the two control methods is subsequently  illustrated through numerical simulations, and various experimental results are presented to demonstrate the effectiveness of the proposed controller to mitigate disturbance amplification and achieve platoon stability.
\end{abstract}
%-------------------------------------------------------------------------------

\section{Introduction}
In recent decades, thanks to developments in automation, such as the emergence of automated vehicles or automated infrastructures, a tremendous revolution has occurred in transportation systems. The transition phase from using only human-driven vehicles (HDVs) to fully   connected and automated vehicles (CAVs) results in new challenges and creates a strong motivation to study the mixed traffic systems, that include both HDVs and CAVs (e.g., see \cite{harfouch2017adaptive, 9029529} and the references therein). In this direction, new opportunities arise to utilize the potential abilities of CAVs to control a transportation network, manage congestion, and promote the efficiency and safety of traffic systems.

More traditional methodologies for controlling traffic flow employ controllers and actuators at fixed locations, among which variable speed limits (VSLs) and ramp metering (RMs) can be mentioned \cite{papamichail2008integrated}. However, the installation of these actuators is not cost-effective and reduces the flexibility of the control system. On the other hand, the advent of CAVs as mobile actuators--so-called Lagrangian actuators--paves the way for applying traffic flow control in a more effective and flexible manner.  For instance, if all the involved vehicles are CAVs, efficient control strategies, such as adaptive cruise control (ACC) and cooperative adaptive cruise control (CACC), can be employed and lead to a desirable system performance \cite{naus2010string,milanes2013cooperative,alam2015experimental, li2017dynamical}.
Nevertheless, in a mixed traffic system, where the penetration rate of CAVs is less than one, new challenges are introduced that require further theoretical and experimental analyses.

In this direction, we consider a mixed traffic system, including one CAV and numerous HDVs, for the two common scenarios of a ring-road and an open-road and reveal how the single CAV is capable of controlling the entire network. In the following, we first present a review of some relevant works in the literature, and subsequently, we discuss the main contributions of this work.

\subsection{Literature Review}

There are a few experimental studies that verify the emergence of stop-and-go waves in  traffic flow systems. For instance, in the study in \cite{sugiyama2008traffic}, the outcome of a practical experiment on a single-lane ring-road demonstrated that a platoon consisting only of HDVs has the potential to initiate stop-and-go waves. These waves that travel upstream along the road make a uniform flow unstable and create a so-called \emph{phantom traffic jam}. This phenomenon of instability has been studied in the literature from macroscopic \cite{flynn2009self}, cellular automaton \cite{nagel1992cellular}, and microscopic \cite{bando1995dynamical} point of view. The emerging nonlinear waves can be amplified by some effects such as stochastic behaviour of human drivers, lane changing, road characteristics, and ramps, to name a few. In \cite{stern2018dissipation}, a field experiment was conducted to show that utilizing a single CAV in a platoon on a circular roadway can dissipate the undesired waves. Moreover, in \cite{cui2017stabilizing}, through some theoretical analysis, the capability of a single CAV to control the traffic flow on a ring-road was investigated.

In fact, since the platoon is connected, and neighboring vehicles can interact, a sparse number of CAVs--that act as mobile actuators--can  influence  the whole network and stabilize the traffic system. The notions of string and ring stability have been employed here for the stability analysis of interconnected vehicles on a string and on a ring roadway, respectively \cite{swaroop1996string, ploeg2013lp,giammarino2020traffic}. In the same direction, and for mixed traffic systems, the string stability of a mixed platoon of infinite length has been analyzed in \cite{talebpour2016influence}. Furthermore, in  \cite{xie2018heterogeneous}, a linear stability condition has been stated in terms of the penetration rate and spatial distribution of CAVs.

In order to dissipate stop-and-go waves in a mixed traffic system, an appropriate control strategy can be provided that is applied to the CAVs as the controllers. To establish a theoretical analysis for these systems, we should first derive a mathematical model that represents the dynamical behaviour of HDVs. In this direction, there are different car-following models, among which the optimal-velocity-follow-the-leader (OV-FTL) model  \cite{bando1995dynamical} and the intelligent driver model (IDM) \cite{treiber2000congested} can be indicated. These models are nonlinear in principle, but  most of studies in the field of mixed traffic systems utilize a linearized version of the nonlinear dynamics around the equilibrium flow (see e.g., \cite{wu2018stabilizing,delle2019feedback, giammarino2020traffic,zheng2020smoothing, wang2020controllability}).

A fundamental network property that should be checked before designing a controller is its controllability or stabilizability \cite{mousavi2018structural, mousavi2020strong, mousavi2020laplacian}. If a linear system is controllable, then a control signal can be designed to steer the system from any initial state to any final state within finite time. A weaker condition that is necessary for the existence of a controller is the stabilizability of a system. A linear system is stabilizable if with a suitable choice of control signals, all the states remain bounded or converge to constant values \cite{callier2012linear}. Recently, some works in the literature have focused on providing a rigorous controllability and stabilizability analysis for ring-road mixed traffic systems with one single CAV.

In \cite{cui2017stabilizing, zheng2020smoothing}, it is assumed that all HDVs in the platoon are homogeneous, which is a strong assumption for practical scenarios. In fact, the problem where heterogeneity of HDVs is considered is closer to reality, but theoretically more challenging. In \cite{wang2019controllability, wang2020controllability}, a controllability analysis for mixed traffic systems, including one CAV and a number of heterogeneous HDVs, is provided. It is stated that under a restrictive condition on  parameters of the dynamic model, all nonzero eigenvalues of the system are controllable; while there exists only one uncontrollable eigenvalue at origin. More recently, in \cite{chen2021mixed}, by considering the similar condition on the parameters of the system, the controllability of a so-called ``$1+n$'' mixed platoon, forming a string at a signalized intersection, is provided. Moreover, in \cite{wang2020leading}, by defining a new notion of leading cruise control, the controllability of a platoon along a string is analyzed. None of  these works provide the controllability analysis of a mixed platoon  in the most general case and without assuming any  constraint on the system parameters.  %However, through numerical simulations, it can be verified that this statement is inaccurate. In fact, for different platoon sizes, one can demonstrate that even if the above-mentioned condition holds, the system has more than one uncontrollable eigenvalues which are nonzero.

%In this paper, we consider a mixed traffic system with one CAV and numerous heterogeneous HDVs. The mixed platoon is studided in two cases of a ring-road and an open-road.  The traffic system is modeled at a microscopic scale, and the dynamics of HDVs are represented by the optimal velocity model. Furthermore, the acceleration of the CAV is directly controlled. For this system, we first establish a stabilizability analysis for the two cases of a ring-road and an open-road in Section III. Moreover, since the goal is to synthesize an output-feedback controller, the detectability analysis of the system is also necessary, which is presented subsequently. Note that

There are numerous works in the literature, proposing various control strategies to stabilize a mixed platoon on an open-road or a ring-road \cite{jin2016optimal,wu2017emergent,cui2017stabilizing, gong2018cooperative, stern2018dissipation,kreidieh2018dissipating,  di2019cooperative} . However, in most of these works, the communication topology of the network which represents the capability of the CAV to receive information from its neighboring vehicles has been neglected \cite{jin2016optimal, vinitsky2018lagrangian,wu2018stabilizing,monteil2018mathcal, di2019cooperative}, and it has been assumed that it can be connected to any vehicle in the platoon. More recently, in \cite{wang2019controllability}, the issue of limited communication has been considered, and a \emph{structured} optimal control has been proposed, which is in general computationally intractable \cite{jovanovic2016controller} and results in a sub-optimal solution. In addition, the HDVs usually do not follow deterministic dynamic models, and there exists uncertainties in the model parameters of HDVs, that can affect the efficiency of the control strategies. Accordingly, it is needed to provide robust control methods that dissipate the perturbations of the traffic flow in the presence of system uncertainties (see e.g., \cite{jin2018experimental, chen2018robust, feng2021robust}).

\subsection{Contributions}

In this paper, we consider a mixed traffic system with one CAV and numerous heterogeneous HDVs. The mixed platoon is studided in two cases of a ring-road and an open-road.  The traffic system is modeled at a microscopic scale, and the dynamics of HDVs are represented by the optimal velocity model. Furthermore, the acceleration of the CAV is directly controlled. For this system, we first establish a stabilizability analysis for the two cases of a ring-road and an open-road. Moreover, since the goal is to synthesize an output-feedback controller, the detectability analysis of the system is also necessary, which is presented subsequently.  %Note that there are numerous works in the literature that propose different control strategies to stabilize a platoon on a ring-road. However, in most of these works, the communication topology of the network which represents the capability of the CAV to receive information from its neighboring vehicles has been neglected \cite{jin2016optimal, vinitsky2018lagrangian,wu2018stabilizing, di2019cooperative}, and it has been assumed that it can be connected to any vehicle in the platoon. More recently, in \cite{wang2019controllability}, the issue of limited communication has been considered, and a \emph{structured} optimal control has been proposed, which is in general computationally intractable \cite{jovanovic2016controller} and results in a sub-optimal solution.
As a real-world scenario,  we also consider the limited communication capability of the CAV in this work. In order to deal with the topological communication constraints, we propose an output-feedback controller that employs the information of a sparse set of vehicles in the control signal design.  Our proposed method offers a robust $H_{\infty}$ controller, that not only increase the efficiency of the CAV, but it also dampens the disturbances appearing as nonlinear waves and improves the performance in the behavior of the entire traffic network.

In summary, the main contributions of this work are listed as follows:

\begin{itemize}
    \item We analyze the controllability and observability of a mixed traffic system with one CAV and numerous heterogeneous HDVs for both common scenarios of a ring-road and an open-road in the most general case. In fact, unlike the existing works in the literature, such as \cite{wang2020controllability}, that investigates the stabilizability under restrictive parameter constraints,  we prove that, for any value of system parameters, the mixed platoon is stabilizable. This analysis verifies the ability of the single CAV to make the states of all HDVs converge to desired values. Further, since we aim to  synthesize a dynamic output controller  that can utilize the states of only a subset of HDVs, we prove in this work that the mixed traffic system is also detectable. %since the system is detectable, synthesizing an output-feedback controller that can only use the states of a subset of HDVs is feasible.
    \item In order to dampen the undesired perturbations occurring in a mixed traffic flow where there is no uncertainties in the system parameters, we propose a  solution based on synthesizing an $H_{\infty}$ output dynamic feedback controller. This controller also tackles the issue of the CAV's communication constraints. To the best of our knowledge, this is the first time in the literature that an output  dynamic controller is utilized for the control of a mixed platoon. Unlike some existing control strategies, e.g.,\ \cite{li2017dynamical, orosz2016connected}, that aim to increase the local efficiency around the CAVs, our method offers a  controller that improves the performance in the behavior of the entire traffic network. More importantly, as we  consider a higher degree of freedom in designing the control strategy, our proposed control method leads to a global optimal solution, while the structural control method in \cite{wang2020controllability} results is a sub-optimal one.
    \item As the next and the more practical scenario, we consider the model mismatch and parametric uncertainties in the dynamics of HDVs and provide a robust control strategy that can smoothen the traffic flow in the presence of disturbances and uncertainties.
\end{itemize}

%Therefore, the main contributions of this work compared to \cite{wang2020controllability} can be summarised as: 1) Rather than the structural optimal control method, presented in \cite{wang2020controllability}, we propose a new solution based on synthesizing an output dynamic feedback controller to tackle the problem of communication constraints. To our best knowledge, this is the first time in the literature that an output  dynamic controller is utilized for the control of a platoon. This approach in fact leads to an optimal solution, while structural control method in \cite{wang2020controllability} results in a sub-optimal one. 2) Unlike \cite{wang2020controllability} that proves the stabilizability of a heterogeneous platoon under restrictive conditions on the value of system parameters, we demonstrate the system stabilizability for a general system with any arbitrary parameters. 3) In order to utilize a dynamic output controller, we prove in this work that the system is also detectable.

     \subsection{Outline}
     The rest of this paper is organized as follow. In Section II, the model of a mixed traffic system is presented, and the main problems of this work are formulated. Section III discusses the stabilizability and detectability of a mixed platoon with a single CAV for both ring-road and open-road setups. In section IV, we propose an output dynamic controller for a mixed traffic flow that has no uncertainty in the system parameters. Section V establishes a control strategy that smoothen the traffic flow in the presence of uncertainties. In Section VI, numerical validations of the results as well as a numerical comparison with some of existing works in the literature are provided. Finally, Section VII concludes the paper.

\section{Preliminaries}

In this section, we first present a dynamic model for a mixed traffic system on a single-lane ring-road, and then we formulate the problem.

We denote the set of real and complex numbers by $\mathbb{R}$ and $\mathbb{C}$, respectively. For $a\in \mathbb{C}$, $\mathrm{Re}(a)$ represents its real part.
%For a vector $v$, $v_i$ is its $i$th entry; for a matrix $M$,  $M_{ij}$ is the entry in row $i$ and column $j$. A subvector $v_X$ is comprised of $v_i$, for $i\in X$, ordered lexiographically.
We denote the transpose of matrix $M$ by $M^T$. Also, $\mathrm{det}(M)$ represents its determinant. For a vector space $\mathcal{V}$, $\mathrm{dim}(\mathcal{V})$ indicates its dimension.
 %$\mathbf{1}_n$ denotes the vector of all ones in $\mathbb{R}^n$.
%
The identity matrix is denoted by $I$, and  its $j$-th column is designated by $e_j$. Also, $0_{n\times m}$ represents an $n\times m$ zero matrix. We also show by $0$ a zero matrix of an appropriate dimension. For $\alpha_1, \ldots, \alpha_n\in\mathbb{R}$, $A=\mathrm{diag}(\alpha_1,\ldots,\alpha_n )$ is an $n\times n$ diagonal matrix whose diagonal elements are $\alpha_1, \ldots, \alpha_n$. For a matrix $M\in\mathbb{R}^{n\times n}$, where $M=M^T$,  $M\succ 0$ (resp., $M\succeq 0$) implies that $M$ is a positive definite (resp., positive semi-definite) matrix. $||.||$ denotes the Euclidean norm
of vectors and the induced norm of matrices. Also, $||.||_{\mathcal{F}}$ denotes Frobenius norm of matrices.

\begin{fact}{\cite{horn2012matrix}} For any matrix $A \in \mathbb{R}^{n\times n}$, we have $||A||<||A||_{\mathcal{F}}$, where $||A||_{\mathcal{F}}=(\sum_{i=1}^{n}\sum_{j=1}^{n} a_{ij}^2)^{\frac{1}{2}}$. \label{fact}\end{fact}

\begin{fact}{\cite{horn2012matrix}} For any matrix $A \in \mathbb{R}^{n\times n}$, $||A||\leq 1$ if and only if $A^T A \preceq I$. \label{fact1}\end{fact}
\subsection{Modeling a Mixed Traffic System}

We study a mixed traffic system that is a network of $n$ vehicles, including one CAV and $n-1$ HDVs. We consider two cases of a single-lane ring-road and a single-lane open-road with length $D$. In Fig.~\ref{net}(a), a schematic diagram of this network is illustrated, where the red car denotes the CAV, indexed with $1$, and all others are HDVs. The position and  velocity of vehicle $i$ are denoted by $p_i$ and $v_i$,  respectively. We define as $s_i=p_{i-1}-p_i$ the back-to-back distance of the $i$-th vehicle from the $i-1$-th vehicle. %.the spacing between vehicles $i$  and $i-1$.

There are different models in the literature to represent the car-following dynamics of human-driven vehicles (see e.g., \cite{treiber2000congested, delle2019feedback,bando1995dynamical}). For instance, the optimal velocity model (OVM) \cite{bando1995dynamical} can be described as:
\begin{equation}
\begin{aligned}
\dot{s}_i(t)&=v_{i-1}(t)-v_i(t)\\
\dot{v}_i(t)&=H_i( v_i(t), s_i(t), \dot{s}_i(t)),
\label{eq1}
\end{aligned}
\end{equation}
where $H_i(\cdot)$ is the acceleration of vehicle $i$, that is a nonlinear function of its velocity $v_i$, the spacing $s_i$, and the relative velocity $\dot{s}_i$. Note that unlike most of the works in literature (e.g., \cite{cui2017stabilizing, jin2016optimal, 9029529}), in this paper, HDVs are assumed to be heterogeneous, and thus, the dynamics of each vehicle $i$ is described by a distinct nonlinear function $H_i(\cdot)$. One can see that at the equilibrium point of dynamics (\ref{eq1}), all vehicles have the same velocity $v^*$. Moreover, since we have $\dot{v}^*=0$, the spacing $s^*_i$ is computed by $0=H_i(v^*, s^*_i, 0)$.
%%%%%%%%%%%%%%%%%%%%%%%%%%%%%%%%%%555
%%%%%%%%%%%%%%%%%%%%%%%%%%%%%%%%%
\begin{figure}[t]
\includegraphics[width=.5\textwidth]{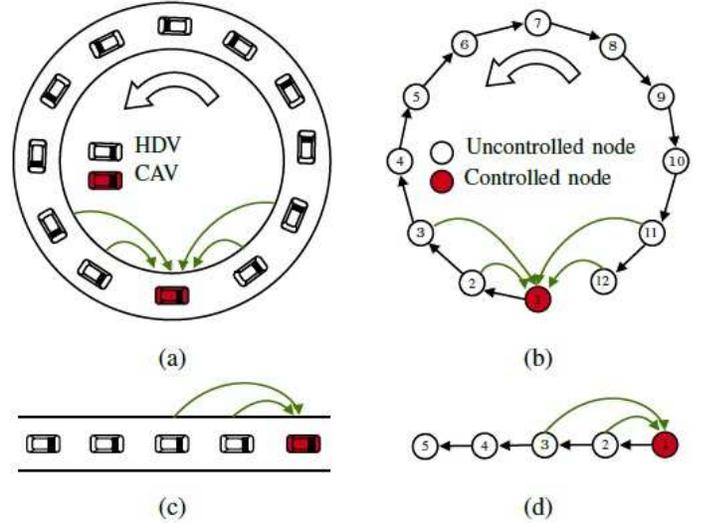}
\centering
\caption{(a)--(b) Schematic of a mixed traffic system on a ring-road and the corresponding graph; (c)--(d) Schematic of a mixed traffic system on an open-road and the corresponding graph. }
\label{net}
\end{figure}
Now, let us define the state error $x_i^T=\begin{bmatrix}
\bar{s}_i & \bar{v}_i
\end{bmatrix}=\begin{bmatrix}
s_i-s^*_i & v_i-v^*
\end{bmatrix}$. Then, by linearization of (\ref{eq1}) around the equilibrium point $\begin{bmatrix}
s^*_i & v^*
\end{bmatrix}^T$, for $i=2, \ldots, n$, one can derive a linear time-invariant (LTI) model for the $i$-th HDV as:
\begin{equation*}
\begin{aligned}
\dot{\bar{s}}_i(t) &= \bar{v}_{i-1}(t)-\bar{v}_i(t)\\
\dot{\bar{v}}_{i}(t) &=\beta_{i1} \bar{s}_i(t) -\beta_{i2} \bar{v}_i(t) +\beta_{i3} \bar{v}_{i-1}
\end{aligned},
\end{equation*}
where
\begin{equation*}
\beta_{i1}=\dfrac{\partial{H_i}}{\partial{s_i}}, \ \ \beta_{i2}=\dfrac{\partial{H_i}}{\partial{\dot{s}_i}}-\dfrac{\partial{H_i}}{\partial{v_i}}, \ \ \beta_{i3}=\dfrac{\partial{H_i}}{\partial{\dot{s}_i}},
\end{equation*}
computed at the equilibrium point. Due to some physical constraints imposed by the behavior of HDVs in practice \cite{cui2017stabilizing}, one should consider
\begin{equation} \beta_{i1}>0, \:\:\: \beta_{i2}>0, \:\:\:  \beta_{i3}>0. \label{beta}\end{equation}

Now, for the dynamics of the single CAV, we consider two cases of a ring-road and an open-road:
1) Corresponding to the case of a ring-road, the dynamics of the single CAV whose acceleration can be directly controlled are given by:
\begin{equation}
\begin{aligned}
\dot{\bar{s}}_1(t) &= \bar{v}_{n}(t)-\bar{v}_1(t)\\
\dot{\bar{v}}_{1}(t) &=u(t)
\end{aligned},
\end{equation}
where $u(t)\in \mathbb{R}$ is the control signal. 2) In the case of an open road, we have:
\begin{equation}
\begin{aligned}
\dot{\bar{s}}_1(t) &= u_1(t)-\bar{v}_1(t)\\
\dot{\bar{v}}_{1}(t) &=u_2(t)
\end{aligned},
\end{equation}
where $u_1(t), u_2(t)\in\mathbb{R}$ are two external inputs, and $u(t)=\begin{bmatrix}
u_1 & u_2\end{bmatrix}^T\in \mathbb{R}^{2}$. In fact, in the case of an open-road, both acceleration and velocity of the CAV are directed by external inputs.

Now, by defining the aggregated vector of states of all vehicles as $x=\begin{bmatrix}
x_1^T & x_2^T & \ldots & x_n^T
\end{bmatrix}^T\in \mathbb{R}^{2n}$, one can derive the following LTI dynamics for the overall system:
\begin{equation}
\dot{x}(t)=Ax(t)+Bu(t),
\label{main}
\end{equation}
\begin{equation}A=\begin{bmatrix}
J_1 & 0 & \ldots & \ldots & 0 & J_2\\
A_{21} & A_{22} & 0 & \ldots & \ldots & 0\\
0 & A_{31} & A_{32} & 0 & \ldots & 0\\
\vdots  & \ddots & \ddots & \ddots & \ddots & \vdots\\
0 & \ldots & \ldots & 0 & A_{n1} & A_{n2}
\end{bmatrix}, \:\:\: B=\begin{bmatrix}
B_1\\ B_2\\ B_2\\ \vdots \\ B_2
\end{bmatrix}, \label{A}\end{equation}
with $J_1=\begin{bmatrix}
0 & -1\\
0 & 0
\end{bmatrix}$, and for $i=2,\ldots, n$, we have:
\begin{equation}A_{i1}=\begin{bmatrix}
0 & 1\\ 0 & \beta_{i3}
\end{bmatrix}, \:\:\: A_{i2}=\begin{bmatrix}
0 & -1 \\ \beta_{i1} & -\beta_{i2}
\end{bmatrix}.\label{Ai}\end{equation}
Now, for the case of a ring-road, we have $B\in\mathbb{R}^{n\times 1}$ and  define:
\begin{equation}
 J_2=\begin{bmatrix}
0 & 1\\0 & 0
\end{bmatrix}, \:\:\:B_1=\begin{bmatrix}
0\\1
\end{bmatrix}, \:\:\: B_2=\begin{bmatrix}
0\\0
\end{bmatrix}\label{ringR}.\end{equation} Moreover, for  an open-road, we have $B\in\mathbb{R}^{n\times 2}$, and $J_2$, $B_1$,  and $B_2$ are defined as:
\begin{equation}
 J_2=\begin{bmatrix}
0 & 0\\0 & 0
\end{bmatrix}, \:\:\:B_1=\begin{bmatrix}
1& 0\\0 & 1
\end{bmatrix}, \:\:\: B_2=\begin{bmatrix}
0 & 0\\0 & 0
\end{bmatrix}\label{openR}.\end{equation}

Finally, we note that the CAV can receive the state information associated with only a number of HDVs, due to its communications constraints. For example, in Fig.~\ref{net}, the green links show HDVs whose information is available to the CAV. In this direction, an output vector $y(t)\in \mathbb{R}^{2m}$ is defined, that includes the states of the CAV together with the states of all HDVs whose information is accessible by the CAV. For $k=1, \ldots, m$, let $j_k$ be the index of the vehicle whose state information, i.e.\ $\bar{s}_{j_k}$ and $\bar{v}_{j_k}$, can be directly measured and observed. Consequently, we have:
\begin{equation}
y(t)=Cx(t),
\label{main11}
\end{equation}
where
$C^T=\begin{bmatrix}
 e_{(2*j_1-1)} & e_{2*j_1} &\ldots & e_{(2*j_m-1)} & e_{2*j_m}
\end{bmatrix}$. Note that we can have $j_1=1$, since the CAV has access to its own information.

In this paper, we design a dynamical output-feedback controller that has the following dynamics:
\begin{equation}
\begin{aligned}
\dot{x}_k(t)&=A_kx_k(t)+B_ky(t)\\
u(t)&=C_kx_k(t)
\end{aligned},\label{cont}
\end{equation}
where $A_k \in \mathbb{R}^{2n\times 2n}$, $B_k \in \mathbb{R}^{2n\times 2m}$ and $C_k \in \mathbb{R}^{1\times 2m}$, and $x_k \in \mathbb{R}^{2n}$ is state of the controller.

\textbf{Problem 1:} In order to ensure the existence of an  output-feedback controller, based on the separation principle \cite{callier2012linear}, the first goal of this paper is to  prove the stabilizability of pair $(A,B)$ and detectability of pair $(A,C)$ in (\ref{main}) and (\ref{main11}).

Pair $(A,B)$ is stabilizable if the uncontrollable modes are all stable. Similarly, pair $(A,C)$ is detectable if the unobservable modes are all stable. We can use the Popov-Belevitch-Hautus (PBH) test for checking the controllability and observability of a specific eigenvalue.

 \color{black}{\begin{pro}[\cite{sontag2013mathematical}]
An eigenvalue $\lambda$ of $A$ for a pair $(A,B)$ (resp., $(A,C)$) is controllable (resp., observable)
if and only if for all nonzero $\rho$ for which $\rho^{T}A=\lambda \rho^T$ (resp., $A\rho=\lambda \rho$), $\rho^{T}B\neq0$ (resp., $C\rho \neq 0$).
\label{p}
\end{pro}}

\textbf{Problem 2:} The next problem of this work is dedicated to computing  matrices $A_k$, $B_k$, $C_k$, and $D_k$ in (\ref{cont}), for the mixed traffic system described in (\ref{main}) and (\ref{main11}),  such that the undesired  perturbations in the traffic flow are dissipated. We first solve this problem  without considering any uncertainty in the parameters of heterogeneous vehicles. Subsequently, we assume that there might be uncertainties in the dynamical model of HDVs, and we compute a robust output-feedback controller  in the presence of disturbance and parametric uncertainty.

Accordingly, in the next section, we  analyze the stabilizability and detectability of the mixed traffic system (\ref{main}) and (\ref{main11}), and subsequently, we design a dynamic controller.

\section{Stabilizability and Detectability}\label{st_det}
In this section, the stabilizability and detectability of dynamical system (\ref{main}) and (\ref{main11}) are discussed.

\subsection{Stabilizability Analysis}

In order to prove the stabilizability of the mixed traffic system, we first consider a ring-road and  prove that there is only one uncontrollable eigenvalue at the origin.
\begin{pro}
In the case of  a ring-road, the pair $(A,B)$  in (\ref{main}), where $B$ is defined in (\ref{ringR}), has only one uncontrollable eigenvalue at the origin.
\label{pp}
\end{pro}

\emph{Proof:}  Let $\rho=\begin{bmatrix}
\rho_1^T & \ldots &\rho_n^T
\end{bmatrix}\in \mathbb{R}^{2n}$, where $\rho_i^T=\begin{bmatrix}
\rho_{i1} & \rho_{i2}
\end{bmatrix}^T$. Considering the expression of $A$ in (\ref{A}), equation $\rho^TA=0$ leads to
\begin{equation}
\begin{aligned}
&\rho^T_iA_{i2}+\rho^T_{i+1}A_{(i+1)1}=0, \:\:\: i=2,\ldots, n-1\\
& \rho_1^TJ_1+ \rho_2^TA_{21}=0\\
& \rho_1^TJ_2+\rho_n^TA_{n2}=0
\label{eqq}
\end{aligned}.
\end{equation}
Now, from (\ref{eqq}), one can derive
\begin{equation*}
\begin{aligned}
&\rho_{i2}=0, \:\:\:\:\:\:\:\: \:\:\:  \:\:\:  \:\:i=2,\ldots,n \\
& \rho_{i1}=\rho_{(i+1)1}, \:\:\: \:\: i=1,\ldots, n-1.
%\label{eqq}
\end{aligned}
\end{equation*}
Thus, a left eigenvector $\rho$ of $A$ associated with $\lambda=0$ can be written as $\rho=[\alpha, \beta, \alpha,0,\ldots, \alpha,0]^T$. Now, assume that $\lambda=0$ is an uncontrollable eigenvalue. Thus, from Proposition 1, one can see that $\rho^TB=0$, which essentially means $\rho_{12}=0$. Now, define $V_u=\{\rho\in \mathbb{R}^{2n}: \rho^TA=0, \rho^TB=0\}$. Then, one can see that $V_u=\{\rho\in \mathbb{R}^{2n}: \rho=[\alpha, 0, \alpha,0,\ldots, \alpha,0]^T, \mathrm{for\:\: some} \:\;\alpha \neq 0\}$. Thus, $\mathrm{dim}(V_u)=1$, which implies that there is only one uncontrollable eigenvalue of $A$ at the origin.
\carre

Now, consider  an open-road. Next, we show that in this case, there is no uncontrollable eigenvalue at the origin.

\begin{pro}
In the case of  an open-road, the pair $(A,B)$  in (\ref{main}), where $B$ is defined in (\ref{openR}), has no uncontrollable eigenvalue at the origin.
\label{ppadd}
\end{pro}

\emph{Proof:} In this case,  from (\ref{eqq}),  we have:
\begin{equation*}
\begin{aligned}
&\rho_{i2}=0, \:\:\:\:\:\:\:\: \:\:\:  \:\:\:  \:\:i=2,\ldots,n \\
& \rho_{i1}=\rho_{(i+1)1}, \:\:\: \:\: i=1,\ldots, n-1\\
& \rho_{n1}=0.
%\label{eqq}
\end{aligned}
\end{equation*}
Therefore, a nonzero left eigenvector $\rho$ of $A$ associated with $\lambda=0$ has a form as $\rho=\begin{bmatrix}0 & \beta & 0 &\ldots &0\end{bmatrix}^T$. Assume that we have an uncontrollable eigenvalue at the origin. Thus, from Proposition 1, $\rho^TB=0$ leads to $\rho_{11}=\rho_{12}=0$, implying that $\rho=0$, that is a contradiction.
\carre

Next, it suffices to show that  all unstable eigenvalues  of $A$ in (\ref{main}), if any, %which lies in the right half-plane %, i.e. any $\lambda$ with $\mathrm{Re}(\lambda)>0$,
  are controllable.

  \begin{lem}
  Let $\lambda\in \mathbb{C}$ with $\mathrm{Re}(\lambda)>0$, and for $i=2, \ldots, n$, let $A_{i2}$ be defined as (\ref{Ai}). Then, the matrix $(\lambda I-A_{i2})$ is nonsingular.
  \label{lem1}
  \end{lem}

\emph{Proof:} For $i=2,\ldots, n$, define $\mathcal{S}_i(\lambda)=\mathrm{det}(\lambda I -A_{i2})=\lambda^2+\beta_{i2}\lambda +\beta_{i1}$.  Based on this equation, the sum of eigenvalues of $A_{i2}$ equals $-\beta_{i2}$, and their product is  $\beta_{i1}$. If the roots of $\mathcal{S}_i(\lambda)$ are real, then,  since $\beta_{i1},\beta_{i2}>0$ from (\ref{beta}),  both the roots are negative. Moreover, if the roots are complex and are written as $a+ib$ and $a-ib$, then since $2a=-\beta_{i2}$, we have $a<0$.   Therefore, because we have $\mathrm{Re}(\lambda)>0$,  one can conclude that $\mathcal{S}_i(\lambda)\neq 0 $, implying that for $i=2,\ldots, n$, $(\lambda I-A_{i2})$ is nonsingular.
\carre

\begin{theo}
  For both cases of a ring-road and an open-road, the pair $(A,B)$  associated with  the mixed traffic system described in (\ref{main}) is stabilizable.
  \label{T}
\end{theo}

\emph{Proof:} Let $\lambda$ be an eigenvalue of $A$ with $\mathrm{Re}(\lambda)>0$, and let also $\rho=\begin{bmatrix}
\rho_1^T &\ldots &\rho_n^T
\end{bmatrix}^T\in \mathbb{R}^{2n}$ be its nonzero left eigenvector, where $\rho_i^T=\begin{bmatrix}
\rho_{i1}&\rho_{i2}
\end{bmatrix}$. Define $A_{11}=J_2$ and $A_{12}=J_1$. Then, for $i=1,\ldots, n$, equation $\rho^TA= \lambda\rho^T$ implies that
\begin{equation}
\rho_i^T(\lambda I -A_{i2})=\rho_{i+1}^T A_{(i+1)1}.
\label{11}
\end{equation}
Now, one can see that $\mathrm{det}(\lambda I -A_{12})=\lambda^2$, which is nonzero, since $\lambda\neq 0$. Moreover, since $\mathrm{Re}(\lambda)>0$, Lemma 1 implies that $(\lambda I-A_{i2})$ is nonsingular. Accordingly, for $i=1, \ldots, n$, one can rewrite (\ref{11}) as
\begin{equation}
\rho_i^T=\rho_{i+1}^T A_{(i+1)1}(\lambda I -A_{i2})^{-1}.
\label{1212}
\end{equation}

Let $L_i=(\lambda I -A_{i2})^{-1}$. Now, by recursively employing equation (\ref{1212}) for $i=1,\ldots, n$, we can derive
\begin{equation}
\begin{aligned}
\rho_1^T&= \rho_{2}^T A_{21} L_1 \\&= \rho_{3}^T A_{31} L_2 A_{21} L_1 \\&= \ldots\\ &= \rho^T_n A_{n1} L_{n-1} A_{(n-1)1} \ldots L_2 A_{21} L_1 \\&= \rho^T_1 A_{11} (L_n A_{n1}) \ldots (L_2 A_{21}) L_1.
\end{aligned}
\label{100}
\end{equation}
In the case an open-road, since $A_{11}=0$, one can conclude from (\ref{100}) that $\rho_1=0$. Now, we want to prove that $\rho=0$ holds for a ring-road as well.
Note that we have
\begin{equation}
L_i A_{i1}= \dfrac{1}{s^i_1}\begin{bmatrix}
0 & s^i_2\\
0 & s^i_3
\end{bmatrix}, \:\:\:\:i=2,\ldots, n,
\label{101}
\end{equation}
where $s^i_1=\lambda^2+\beta_{i2}\lambda +\beta_{i1}$, $s^i_2=\lambda+\beta_{i2}-\beta_{i3}$, and $s^i_3=\beta_{i3}\lambda+\beta_{i1}$. As shown before, $s_1^i\neq 0$. %Moreover,  since $\beta_{i1},\beta_{i3}>0$ from (\ref{beta}) and $\mathrm{Re}(\lambda)>0$, we have $s^i_3\neq 0$.
Now, we substitute (\ref{101}) into (\ref{100}), and for the case of a ring-road, we  obtain
\begin{equation}
\begin{bmatrix}
\rho_{11} & \rho_{12}
\end{bmatrix}=\dfrac{\Pi_{i=2}^n s^i_3}{\lambda \Pi_{i=2}^n s^i_1}\begin{bmatrix}
\rho_{11} & \rho_{12}
\end{bmatrix} \begin{bmatrix}
0 & 1\\
 0 & 0
\end{bmatrix}. \label{hel}
\end{equation}
Equation (\ref{hel}) leads to $\rho_{11}=0$ and $\rho_{12}= \frac{\Pi_{i=2}^n s^i_3}{\lambda \Pi_{i=2}s^i_1} \rho_{11}=0$. Thus, $\rho_1=0$. Now, for both cases of a ring-road and an open-road,  by recursively applying (\ref{1212}) for $i=n, n-1, \ldots, 2$, one can conclude that $\rho_i=0$. Therefore, we have $\rho=0$ that contradicts the assumption. Thus, $A$ has no eigenvalue which lies on the right half-plane. In addition, from Proposition 2, for the case of a ring-road, there is only one uncontrollable eigenvalue at origin. Furthermore, Proposition \ref{ppadd} implies that in the case of an open-road, there is no uncontrollable eigenvalue at origin. Hence, in both cases,  the system is stabilizable.
\carre

\begin{rem}
In \cite{wang2020controllability}, it has been stated that the mixed traffic system described in (\ref{main}) is stabilizable if for all $i,k\in\{1,2,\ldots,n\}$, we have $\beta_{k1}^2-\beta_{i2}\beta_{k1}\beta_{k3}+\beta_{i1}\beta_{k3}^2\neq 0$. Moreover, through a different approach, it has been shown that there is only one uncontrollable eigenvalue at origin.  %However, by numerical simulation, it can be verified that even under the above-mentioned constraint, the system has nonzero eigenvalues that are not controllable.
However, in Theorem \ref{T}, without assuming any restrictive constraint,  we have demonstrated the stabilizability of a mixed traffic system with heterogeneous HDVs  for both cases of a ring-road and an open-road and proved that the system is stabilizable  even if for some $i,k\in\{1,2,\ldots,n\}$, we have $\beta_{k1}^2-\beta_{i2}\beta_{k1}\beta_{k3}+\beta_{i1}\beta_{k3}^2= 0$.
%For instance, consider a mixed traffic system with one CAV and $10$ homogeneous HDVs, and let $\beta_{i1}=\beta_{i3}=1$ and $\beta_{i2}=3$ for all $i\in\{1,\ldots, 11\}$. Then, one can see that the system has $5$ uncontrollable eigenvalues; however,  the system is stabilizable, since it has no eigenvalue on the right half-plane.
\end{rem}

%\begin{rem}
%In \cite{wang2020controllability}, by considering a certain constraint on the values of $\beta_{i1}$, $\beta_{i2}$, and $\beta_{i3}$, for $i=1,\ldots, n$,  it has been stated that all the nonzero eigenvalues associated with a mixed traffic system on a ring-road are controllable. Moreover, through a different approach, it has been shown that there is only one uncontrollable eigenvalue at origin.  %However, by numerical simulation, it can be verified that even under the above-mentioned constraint, the system has nonzero eigenvalues that are not controllable.
%On the other hand, in this paper, the stabilizability of a mixed traffic system with heterogeneous HDVs for both cases of a ring-road and an open-road  is studied, and the stabilizability is demonstrated  for the most general case and without assuming any constraint.
%For instance, consider a mixed traffic system with one CAV and $10$ homogeneous HDVs, and let $\beta_{i1}=\beta_{i3}=1$ and $\beta_{i2}=3$ for all $i\in\{1,\ldots, 11\}$. Then, one can see that the system has $5$ uncontrollable eigenvalues; however,  the system is stabilizable, since it has no eigenvalue on the right half-plane.
%\end{rem}

\subsection{Detectability Analysis}

Here, the detectability of the mixed traffic system is studied, and we show that by observing only the states of the single CAV, the detectability of the whole system can be ensured.

\begin{pro}
The zero eigenvalue of the mixed traffic system described in (\ref{main}) and (\ref{main11}), in both cases of a ring-road and an open-road,  is observable even if the CAV has access to only its own states.
\label{p3}
\end{pro}

\emph{Proof:}
 Let $C=\begin{bmatrix}
e_1 & e_2
\end{bmatrix}$, and assume that $\lambda=0$ is not observable. Then, from Proposition 1, $A$ has a nonzero right eigenvector $\rho$, where we have both $A\rho=0$ and $C\rho =0$. Let $\rho=\begin{bmatrix}
\rho_1^T&\ldots&\rho_n^T
\end{bmatrix}^T\in \mathbb{R}^{2n}$, where $\rho_i^T=\begin{bmatrix}
\rho_{i1} &\rho_{i2}
\end{bmatrix}$. From equation $A\rho=0$, one can write:
\begin{equation}
\begin{aligned}
&\rho_{i2} - \rho_{(i+1)2}=0, &i=2,\ldots, n\\
&\beta_{i3}\rho_{(i-1)2}+ \beta_{i1}\rho_{i1}-\beta_{i2}\rho_{i2}=0, & i=2,\ldots, n
\label{12}
\end{aligned}.
\end{equation}

Moreover, from $C\rho =0$, one can conclude that $\rho_{11}=\rho_{12}=0$.
Thus, (\ref{12}) leads to $\rho=0$ in both cases of a ring-road and an open-road, which is a contradiction. Therefore, the zero eigenvalue is observable.
\carre

\begin{theo}
In the both cases of a ring-road and an open-road, the mixed traffic system described in (\ref{main}) and (\ref{main11}) is detectable even if only the states of the single CAV are directly observed.
\end{theo}

\emph{Proof:} Let $C=\begin{bmatrix}
e_1 & e_2
\end{bmatrix}$. Assume that the system is not detectable. Then, $A$ has an eigenvalue $\lambda$  on the right half-plane with a nonzero right eigenvector $\rho$ such that $C\rho=0$.  Denote  $\rho=\begin{bmatrix}
\rho_1^T&\ldots&\rho_n^T
\end{bmatrix}^T\in \mathbb{R}^{2n}$, where $\rho_i^T=\begin{bmatrix}
\rho_{i1}&\rho_{i2}
\end{bmatrix}$. Then, we should have $\rho_{11}=\rho_{12}=0$. Define $A_{11}=J_2$ and $A_{12}=J_1$. Now, from equation $A\rho=\lambda \rho$, for $i=1,\ldots, n$, one can write:
\begin{equation*}
(\lambda I-A_{i2})\rho_i=A_{i1}\rho_{i-1}.
\end{equation*}
Since $\mathrm{Re}(\lambda )>0$, from Lemma 1, one can see that, for $i=2,\ldots,n$, $(\lambda I -A_{i2})$ is invertible. Moreover, $\mathrm{det}(\lambda I-A_{12})=\lambda^2\neq 0$. Let $L_i=(\lambda I -A_{i2})^{-1}$. Then,  for $i=1,\ldots, n$, we can write:
\begin{equation}
\rho_i=L_i A_{i1}\rho_{i-1}.
\label{detect}
\end{equation}
Now, since $\rho_1=0$, recursively using equation (\ref{detect}) for $i=2,\ldots, n$ leads to $\rho=0$, which contradicts the assumption. In addition, based on Proposition 4, the zero eigenvalue is observable. Thus, in summary, one can conclude that the system is detectable.
\carre

\section{Controller Synthesis: Without Uncertainties}\label{kon1}
In this section, we aim to design a dynamic output-feedback controller for the mixed traffic system (\ref{main}) to dissipate undesired perturbations.
%
%\subsection{Controller Synthesis: Without Uncertainty}
In the first step, we assume there is no uncertainty in the model of the traffic system and design an output-feedback controller.

\subsection{Disturbances and Performance Outputs}\label{D1}
 In fact, the perturbations may appear due to lane changes or merges or the stochastic behavior of HDVs in ring-roads with no bottlenecks \cite{sugiyama2008traffic,wang2020controllability}. %, that leads to traffic waves and increase the traffic congestions.

Perturbations are modeled as disturbance signals added to the acceleration of each vehicle. Thus, by defining   $d=\begin{bmatrix}
d_1(t) & \ldots & d_n(t)
\end{bmatrix}^T\in \mathbb{R}^n$ as the disturbance vector and the matrix $$B_d=\begin{bmatrix}
b_d & 0 & \ldots & 0\\
0 & b_d & \ldots & \vdots\\
\vdots & \ddots & \ddots & 0 \\
0 & \ldots & 0 & b_d
\end{bmatrix}\in\mathbb{R}^{2n\times n},$$ with $b_d=[0,1]^T$, the dynamics of the mixed traffic system in the presence of disturbances can be written as
\begin{equation}
\begin{aligned}
\dot{x}(t)&= Ax(t)+Bu(t)+B_d d(t)\\
z(t)&= C_z x(t) + D_z u(t)\\
y(t)&=Cx(t)
\end{aligned},
\label{main1}
\end{equation}
where for the case of a ring-road, one can define  the controlled output (performance output) $z(t)\in\mathbb{R}^{2n+1}$ as   \begin{equation}z(t)=\begin{bmatrix}
\gamma_s \bar{s}_1(t) & \gamma_v \bar{v}_1 & \ldots & \gamma_s \bar{s}_n(t) & \gamma_v \bar{v}_n & \gamma_u u
\end{bmatrix}^T.\label{z}\end{equation} The parameters $\gamma_s, \gamma_v, \gamma_u>0$ denote the penalties associated with the spacing error, the velocity deviation, and the control energy, respectively. Thus, $C_z=\begin{bmatrix}
\mathcal{T}^{\frac{1}{2}} & 0_{2n\times 1}
\end{bmatrix}^T$ includes  the weights of the
states, with $\mathcal{T}^{\frac{1}{2}}=\mathrm{diag}(\gamma_s, \gamma_v, \ldots, \gamma_s, \gamma_v)$, and  $D_z=\begin{bmatrix}
0_{1\times 2n}&  \mathcal{Q}^{\frac{1}{2}}
\end{bmatrix}^T$ is the
weight of the control input, with $ \mathcal{Q}^{\frac{1}{2}}=\gamma_u$. In the case of an open-road, we define $z(t)=\begin{bmatrix}
\gamma_s \bar{s}_1(t) & \gamma_v \bar{v}_1 & \ldots & \gamma_s \bar{s}_n(t) & \gamma_v \bar{v}_n & \gamma_{u_1} u_1 & \gamma_{u_2} u_2
\end{bmatrix}^T\in\mathbb{R}^{2n+2}$. In this case, one has $C_z=\begin{bmatrix}
\mathcal{T}^{\frac{1}{2}} & 0_{2n\times 2}
\end{bmatrix}^T$ and $D_z=\begin{bmatrix}
0_{2\times 2n}& \mathcal{Q}^{\frac{1}{2}}
\end{bmatrix}^T$, with $ \mathcal{Q}^{\frac{1}{2}}=\mathrm{diag}(\gamma_{u_1}, \gamma_{u_2})$.

\subsection{Output-feedback controller}
Consider an output-feedback controller with the dynamics described in (\ref{cont}).
We aim to design a dynamic controller  that   stabilizes the system (\ref{main1}) and minimizes the influence of the disturbance $d$ on the performance output $z$.

By applying the controller (\ref{cont}) to (\ref{main1}), the closed-loop system is expressed as
\begin{equation}
\begin{aligned}
\dot{\bar{x}}(t)&=\bar{A}\bar{x}(t)+\bar{B} d(t)\\
z(t) &=\bar{C} \bar{x}(t)
\end{aligned},
\label{cc}
\end{equation}
where $\bar{x}=\begin{bmatrix}
{x}(t)\\ {x}_k(t)
\end{bmatrix}$, $\bar{A}=\begin{bmatrix}
A & B C_k\\B_k C&A_k
\end{bmatrix}$, $\bar{B}=\begin{bmatrix}
B_d\\0
\end{bmatrix}$ and $\bar{C}=\begin{bmatrix}
C_z&D_z C_k
\end{bmatrix}$.

\subsection{$H_\infty $ Control Problem}\label{se}
Let $T_{zd}$ be the closed-loop transfer function from the disturbance $d$ to the performance output $z$.

\emph{Problem:} Find the matrices $A_k$, $B_k$, and $C_k$ for the controller (\ref{cont}) such that the closed-loop system (\ref{cc}) satisfies the
 inequality \begin{equation}
||T_{zd}||_{\infty}=\max_{d(t)\neq 0} \frac{||z(t)||_2}{||d(t)||_2}<\gamma,
\label{Hi}
\end{equation}  and  $\gamma >0$ is minimized.

  Note that $||T_{zd}||_{\infty}$ denotes the $H_{\infty}$ norm of $T_{zd}$, which measures the largest input-output gain for energy or power input signals \cite{scherer1997multiobjective}.

\emph{Problem solution:}  Based on the bounded real lemma (BRL) \cite{boyd1994linear}, $||T_{zd}||_{\infty}$ is smaller than $\gamma$ if and only if there exists a positive definite
matrix $P\succ 0$, and matrices $A_k$, $B_k$, and $C_k$ satisfying
%
%\emph{Problem 1:} Minimize $\gamma$ such that there exist a positive definite
%matrix $P\succ 0$, and matrices $A_k$, $B_k$ and $C_k$
%satisfying
\begin{equation}
\begin{bmatrix}
\bar{A}^T P+P\bar{A}&P\bar{B}&\bar{C}^T\\ \bar{B}^T P&-\gamma^2I&0\\ \bar{C}&0&- I
\end{bmatrix}\prec 0. \label{lm}
\end{equation}
Now, note that the inequality (\ref{lm}) is not a linear matrix inequality (LMI) with respect to the variables $P$, $A_k$, $B_k$, and $C_k$, because it is not linear with respect to these variables.  In order to extract an LMI for  computing the controller parameters $A_k$, $B_k$, and $C_k$, we  apply a method that has been presented  in \cite{scherer1997multiobjective} in details, and  we provide a summarized description of this method in the following.

Let us partition the  matrices $P$ and $P^{-1}$ as the following form: $$
P=\begin{bmatrix}
Y&N\\ N^T& *
\end{bmatrix}, \;\;\;
P^{-1}=\begin{bmatrix}
X&M\\ M^T& *
\end{bmatrix},$$ where $X,Y\in \mathbb{R}^{2n\times 2n}$. Note that since $P\succ 0$, we have $X,Y \succ 0$.
From the equation $PP^{-1}=I$,  one can  deduce that
$NM^T+YX=I$.  Further, one can find that $P=\Lambda_2 \Lambda_1^{-1}$, where $$\Lambda_1= \begin{bmatrix}
X&I\\M^T&0
\end{bmatrix}, \:\:\: \Lambda_2=\begin{bmatrix}
I&Y\\0&N^T
\end{bmatrix}.$$ Therefore, we have
$$\Lambda_1^T P \Lambda_1=\Lambda_1^T \Lambda_2=\begin{bmatrix}
X&I\\I&Y
\end{bmatrix}\succ0. $$
Substituting $P=\Lambda_2 \Lambda_1^{-1}$ in (\ref{lm}), we obtain
\begin{equation}
\Omega_1=\begin{bmatrix}
\bar{A}^T \Lambda_2 \Lambda_1^{-1}+\Lambda_1^{-T}\Lambda_2^T \bar{A}&\Lambda_1^{-T}\Lambda_2^T\bar{B}&\bar{C}^T\\ \bar{B}^T\Lambda_2 \Lambda_1^{-1}&-\gamma^2I&0\\ \bar{C}&0&- I
\end{bmatrix}\prec 0. \label{si1} \end{equation}
Now, let us define a matrix $\Omega_2$ by the following congruent transformation
of $\Omega_1$
\begin{equation}
\begin{aligned}
\Omega_2&=\begin{bmatrix}
\Lambda_1^T&0&0\\0&I&0\\ 0&0& I
\end{bmatrix}\Omega_1\begin{bmatrix}
\Lambda_1&0&0\\0&I&0\\ 0&0& I
\end{bmatrix}\\&=\begin{bmatrix}
\Lambda_1^{T}\bar{A}^T \Lambda_2 +\Lambda_2^T \bar{A}\Lambda_1&\Lambda_2^T\bar{B}&\Lambda_1^T\bar{C}^T\\ \bar{B}^T \Lambda_2 &-\gamma^2I&0\\ \bar{C}\Lambda_1&0&- I
\end{bmatrix}\prec 0. \label{si2}\end{aligned} \end{equation}
Now, substituting the values of $\Lambda_1$, $\Lambda_2$, $\bar{A}$, $\bar{B}$ and $\bar{C}$,  we obtain the following matrix inequality:
\begin{equation}
\begin{aligned}
\Omega_2&=\begin{bmatrix}
\Omega_{11}&\Omega_{12}&\Omega_{13}&\Omega_{14}\\\Omega_{21}&\Omega_{22}&\Omega_{23}&\Omega_{24}\\ \Omega_{31}&\Omega_{32}& \Omega_{33}&\Omega_{34}\\\Omega_{41}&\Omega_{42}& \Omega_{43}&\Omega_{44}
\end{bmatrix}\prec 0, \label{si3}\end{aligned} \end{equation}
where
\begin{equation}
\begin{aligned}
&\Omega_{11}=
AX+XA^T+BC_k M^T+ MC_k^T B^T\\&\Omega_{12}=\Omega^T_{21}=MA_k^TN^T+XC^TB_k^TN^T+MC_k^TB^TY\\&\quad\quad\quad\quad\:\:\:+X A^TY+A\\&\Omega_{22}=A^T Y+YA+NB_kC+C^T B_k^T N^T\\&\Omega_{13}=\Omega^T_{31}=B_d,\; \Omega_{14}=\Omega^T_{41} =XC_z^T+M C_k^TD^T_z, \\&\Omega_{23}=\Omega^T_{32}=YB_d,\; \Omega_{24}=\Omega^T_{42}=C_z^T ,\;\Omega_{33}=-\gamma^2I,\\& \Omega_{34}=\Omega_{43}=0,\;\Omega_{44}=- I.\end{aligned} \nonumber\end{equation}
Finally, by defining a new set of variables as
\begin{equation}
\begin{aligned}
&\hat{A}=NA_kM^T+NB_k CX+YBC_k M^T+YAX\\ & \hat{B}=NB_k \\ &\hat{C}=C_k M^T\\ &\eta=\gamma^2 \label{nvar}\end{aligned},\end{equation}
 we obtain the LMIs described in (\ref{lmi}) with respect to variables $\eta, X, Y, \hat{A}, \hat{B}, \hat{C}$.
\begin{figure*}[!t]
\normalsize
\begin{equation}
\begin{aligned}  &\min_{\eta,X,Y,\hat{A},\hat{B},\hat{C}}  \eta \\
 \text{subject to:}&\\
&\begin{bmatrix}
X&I\\I&Y
\end{bmatrix}\succ0 \\
&\begin{bmatrix}
AX+XA^T+B\hat{C}+\hat{C}^T B^T&\hat{A}^T+A & B_d &  XC_z^T+\hat{C}^TD_z^T\\ * &A^T Y+YA+\hat{B}C+C^T \hat{B}^T&YB_d&C_z^T\\ *&*&- \eta I&0\\ *&*&*&- I
\end{bmatrix}\prec 0
\end{aligned}
\label{lmi}
\end{equation}
%\hrulefill
%\vspace{4pt}
\end{figure*}
Therefore, if the optimization problem (\ref{lmi}) is feasible, then we can find $\eta$, $X$, $Y$, $\hat{A}$, $\hat{B}$, and $\hat{C}$  and solve the matrix equation $NM^T=I-YX$ for the non-singular matrices $M$ and $N$. Moreover, matrices $A_k$, $B_k$, and $C_k$ in the state-space realization of the  output-feedback controller %which satisfies $||T_{zd}||_{\infty}<\gamma$ and   minimizes $\gamma $,
can be derived based on (\ref{nvar}) as follows:
\begin{equation}
\begin{aligned}
&A_k=N^{-1}(\hat{A}-NB_k CX-YBC_k M^T-YAX)M^{-T},\\ & B_k=N^{-1}\hat{B}, \\ &C_k=\hat{C}M^{-T}.\label{coeff}\end{aligned}
\end{equation}\carre

\section{Controller Synthesis: With Uncertainties}\label{kon2}

In this section, we consider  parametric uncertainties in the dynamic model of the HDVs. In fact, we  assume that $\beta_{i1}$, $\beta_{i2}$, and $\beta_{i3}$, appearing in (\ref{Ai}), are unknown   parameters for  the mixed traffic system.  Then, we synthesize an output-feedback controller with dynamics (\ref{cont}) that stabilizes the entire closed-loop system
in the presences of disturbances and parametric uncertainties.

%\emph{Assumption 1:}
After the linearization of the overall dynamics of the system, the uncertainty is assumed to appear in the system matrix $A$ in   (\ref{A}) as
\begin{equation}
    A=A_{\mathcal{N}}+\Delta{A}.\label{delta}
\end{equation}
The matrix $A_{\mathcal{N}}$
 is the  mean-valued matrix that is constant and known.  One can calculate $A_{\mathcal{N}}$
as $
A_{\mathcal{N}}=\begin{bmatrix}\frac{a_{ij,\min}+a_{ij,\max}}{2}
\end{bmatrix}$,
where $a_{ij,\min}$ (resp., $a_{ij,\max}$) is the minimum (resp., maximum) value that an entry of $A$ in its $i$th row and $j$th column can have.
%matrix $A(t)$.
On the other hand, $\Delta{A}$ represents parametric uncertainties, that is assumed to be
 structurally bounded. One can  write $\Delta A$ as
$$ \Delta{A}=LFR,$$
where  the matrices $L\in\mathbb{R}^{2n\times 2n}$ and $R\in\mathbb{R}^{2n\times 2n}$ are known  constant
matrices. Moreover, $F\in\mathbb{R}^{2n\times 2n}$ is an unknown matrix, which satisfies the following condition:
\begin{equation}F^TF\preceq I .\label{conun}\end{equation}

\begin{rem}
(Computing $L$ and $R$) We can choose $L$ and $R$ as  $L=\varrho I$ and $R=\rho I$. Therefore, $\Delta A=LFR=\rho\varrho F$, and based on Fact \ref{fact},  we get $||F||=(\rho\varrho)^{-1}||\Delta A||\leq (\rho\varrho)^{-1}||\Delta A||_{\mathcal{F}}$. On the other hand, $\Delta A=A-A_{\mathcal{N}}=\begin{bmatrix}\tilde{a}_{ij}\end{bmatrix}$, where  $|\tilde{a}_{ij}|\leq \frac{1}{2} (a_{ij,\max}-a_{ij,\min})$. Therefore, $||F||\leq  (2\rho\varrho)^{-1}(\sum_{i=1}^{2n}\sum_{j=1}^{2n}( a_{ij,\max}-a_{ij,\min})^2)^{\frac{1}{2}}$. Now, if one sets $\rho\varrho=\frac{1}{2}(\sum_{i=1}^{2n}\sum_{j=1}^{2n}( a_{ij,\max}-a_{ij,\min})^2)^{\frac{1}{2}}$, then $||F||\leq 1$, and the condition (\ref{conun}) is satisfied.
\end{rem}

In the following, we show how we can  ensure the existence of an output-feedback controller that attenuates the effect of the disturbance on the performance output, while stabilizing the closed-loop system.

\begin{theo} Consider a  system with dynamics (\ref{main1}), which has uncertainties that are modelled as (\ref{delta}). Then, there exists an output-feedback controller that stabilizes the closed-loop system (\ref{cc}) and minimizes $\gamma$ in the inequality (\ref{Hi}) if  the optimization problem (\ref{lmiun})  feasible. %, where $\eta$
% and the inequality Under the output-feedback controller (\ref{cont}) and  in the presence of  uncertainties, if the optimization problem (\ref{lmiun}) is feasible, then the closed-loop system  satisfies the inequality  \begin{equation}
%||T_{zd}||_{\infty}=\max_{d(t)\neq 0} \frac{||z(t)||_2}{||d(t)||_2}<\gamma,
%\label{Hii}
%\end{equation}  and  $\gamma >0$ is minimized. %, where $\eta=\gamma^2$.
\begin{figure*}[!t]
\normalsize
\begin{equation}
\begin{aligned}  &\min_{\eta,\epsilon_1,\epsilon_2,\epsilon_3,X,Y,\hat{A},\hat{B},\hat{C}}  \eta \\
 \text{subject to:}&\\
&\begin{bmatrix}
X&I\\I&Y
\end{bmatrix}\succ 0 \\
&\begin{bmatrix}\Gamma_{11}&\Gamma_{12}\\ * &\Gamma_{22}
\end{bmatrix}\prec0\\
\Gamma_{11}=&\begin{bmatrix}
A_{\mathcal{N}}X+XA_{\mathcal{N}}^T+B\hat{C}+\hat{C}^T B^T+L(\epsilon_1+\epsilon_2)L^T&\hat{A}^T+A_{\mathcal{N}}&B_d&XC_z^T+\hat{C}^TD_z^T\\ * &A_{\mathcal{N}}^T Y+YA_{\mathcal{N}}+\hat{B}C+C^T \hat{B}^T+\epsilon_3R^TR&YB_d&C_z^T\\ *&*&-\eta I&0\\ *&*&*& -I
\end{bmatrix}\\
\Gamma_{12}=&\begin{bmatrix}
XR^T&0&0&YL&0\\0&R^T&YL&0&XR^T\\0&0&0&0&0\\0&0&0&0&0
\end{bmatrix}\\
\Gamma_{22}=&-\mathrm{diag}(\epsilon_1I,\epsilon_2I,\epsilon_3I, I, I)
\end{aligned}
\label{lmiun}
\end{equation}
%\hrulefill
%\vspace{4pt}
\end{figure*}
\label{uncertain}
\end{theo}

\emph{Proof of Theorem \ref{uncertain}:} See Appendix.

\begin{rem}
By solving the optimization problem (\ref{lmiun}),  $\eta$, $\epsilon_1,\epsilon_2,\epsilon_3$, $\hat{A}$, $\hat{B}$, $\hat{C}$, $X$, and  $Y$ can be computed. We note that $\eta=\gamma^2$, where $\gamma$ is an upper bound of the $||T_{zd}||_{\infty}$. In addition, we can solve the matrix equation $NM^T=I-YX$ for the non-singular matrices $M$ and $N$. Then, the
matrices $A_k$, $B_k$, and $C_k$ in the state-space realization of the output-feedback controller can be obtained as
\begin{equation*}
\begin{aligned}
&A_k=N^{-1}(\hat{A}-NB_k CX-YBC_k M^T-YA_{\mathcal{N}}X)M^{-T}\\ & B_k=N^{-1}\hat{B} \\ &C_k=\hat{C}M^{-T}.\end{aligned}
\end{equation*}
\end{rem}

%\begin{rem}

%\end{rem}

\section{Simulation Results}

In this section, the efficiency of the proposed control strategies are validated through numerical simulations. The minimization problems (\ref{lmi}) and (\ref{lmiun}) are solved in
YALMIP interface for MATLAB with   Sedumi solver. For a better comparison of the results, we simulate an experimental setup that is  similar to the ones considered in \cite{wang2019controllability, wang2020controllability}, where a ring-road with a circumference $D=400\:\mathrm{m}$ and $20$ vehicles has been studied (since the results for the case of an open-road is analogous to the case of a ring-road, we illustrate here the simulation results only for a ring-road).  We assume  that the  1\textsuperscript{st} vehicle can be a CAV that has access to the state information of the five HDVs ahead and the five HDVs behind.

\subsection{Simulation setup}
Through using an optimal velocity model (OVM), any nonlinear function $H_i(\cdot)$ in (\ref{eq1}), $i=2,\ldots, 20$, that describes the acceleration function of the $i$th HDV  is written as
%
%\begin{equation*}
$H_i(\cdot)=\alpha_i(\mathcal{V}_i(s_i(t))-v_{i}(t)) + \theta_i \dot{s}_i(t) \label{hii}$, where $\alpha_i$ and $\theta_i$ are sensitivity coefficent in an OVM. Moreover, $\mathcal{V}_i(s_i)$ is the desired speed of HDV $i$ that is a function of spacing $s_i$. Due to heterogeneity of HDVs, we set $\alpha_i=0.6 + U[-0.1, 0.1]$ and $\theta_i=0.9+U[-0.1,0.1]$, where $U[a,b]$ represents a uniform distribution function that take values from the interval $[a,b]$. Moreover, we define $\mathcal{V}_i(s_i)$ as a piecewise function
%\end{equation*}
%
$$
\mathcal{V}_i(s_i)=
\begin{cases}
0, & s_i\leq s_{i,st},\\
h_{i,v}(s_i), & s_{i,st} < s_i < s_{i,go},\\
v_{i,\max}, & s_i\geq s_{i,go},
\end{cases}
$$
where we set $s_{i,st}=5$, $s_{i,go}=35+U[-5,5]$, $v_{i,\max}=30$, and $h_{i,v}(s_i)$ is a nonlinear function chosen as   \cite{jin2016optimal}
$$h_{i,v}(s_i)=\dfrac{v_{i,\max}}{2}(1-\cos(\pi \dfrac{s_i-s_{i,st}}{s_{i,go}-s_{i,st}})).$$

In order to ensure the safety and  prevent collisions, every vehicle is also assumed to be equipped with an automatic braking system, described as
$$ \dot{v}(t)= a_{\min} \; \;\;\;\; \; \text{if} \; \;\;\;\; \;\dfrac{v_{i}^2(t)-v^2_{i-1}(t)}{2s_i(t)}\geq |a_{\min}|,$$ where $a_{\min}=-5\:\mathrm{m}/\mathrm{s}^2$. Moreover, we set the maximum acceleration of any vehicle as $a_{\max}=2 \:\mathrm{m}/\mathrm{s}^2$.

\subsection{Stabilizability Verification}

As the first scenario, assume that all vehicles are randomly distributed along the ring-road and start their movement with initial velocity $v_i(0)$ from the distribution $15+U[-4,4] \:\mathrm{m}/\mathrm{s}$. Notice that we consider no uncertainties in the dynamic model of the system in this case. First, assume that all vehicles of the mixed traffic system are HDVs. In Fig.~\ref{pert11}(a), it can be seen that  multiple  perturbations occur in this system, which are amplified over time and generate an unstable nonlinear wave moving upstream the traffic flow.

Next, assume that the 1\textsuperscript{st}  vehicle is a CAV that is controlled by an output-feedback controller described in Section \ref{se}. First, we set the equilibrrium velocity as $v^*=15 \:\mathrm{m}/\mathrm{s}$. Then, the equilibrium spacing  of vehicle $i$ is obtained by solving the equation $0=H_i(15, s_i^*, 0)$. Note that in a circular path, the sum of spacing of all vehicles should equal to the circumference $D$. In other words, we should have $\sum_{i=1}^2 s_i(t)=\sum s_i^*=400$. Now, one can see in Fig.~\ref{pert11}(b) that the perturbations can be attenuated within a short time. As the next experiments, we change the equilibrium velocity to $16 \:\mathrm{m}/\mathrm{s}$ and $14 \:\mathrm{m}/\mathrm{s}$, respectively. Notice that in this case, the parameters of the linearized system (\ref{main}) that are computed around the equilibrium point will be changed. In these two cases, one can observe in Figs.~\ref{pert11}(c,d) that the single CAV is still capable of stabilizing traffic flow and steering it towards the new equilibrium points, which verifies the stabilizability of the mixed traffic system with a single CAV.

%Robustness of the proposed controller is tested for the  system  against unknown  disturbances,
 %and parametric
%uncertainties.  Details of each study scenarios are discussed in the following sub-sections.
 %\subsection{Robustness against   disturbances}
 %The weights in penalty vectors, are selected such that the LMI in problem (\ref{lmi}) is feasible and maintain small system oscillation. The system oscillation increases when small weights are considered for
%the  states. Also, smaller values result as slow response of
%the controller. To keep the system response fast, the state
%weights  are chosen to be $\gamma_s=0.03$, $\gamma_v=0.15$.
%The control effort weight is kept unity $\gamma_u=1$.

\begin{figure}[t]
\includegraphics[width=.5\textwidth]{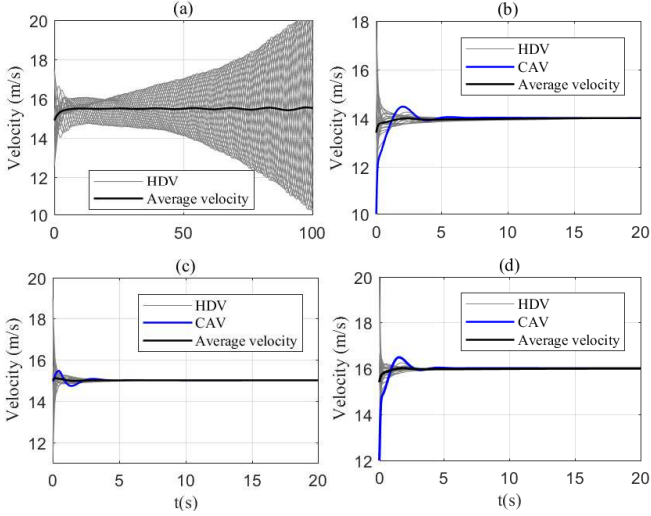}
\centering
\caption{(a) Velocity profile of vehicles when they are all HDVs; (b)--(d) Velocity profile of vehicles when  the 1\textsuperscript{st} vehicle is a CAV under the proposed control in Section \ref{se}, and $v^*$ is 15, 16, and 14 for (b), (c), and (d), respectively.}
\label{pert11}
\end{figure}

\subsection{Robustness Against   Disturbances}
 In this part, we aim to illustrate the performance of the output-feedback controller, proposed in Section \ref{se}, to dampen the undesired disturbances.

 We note that the parameters $\gamma_s$, $\gamma_v$, and $\gamma_u$ in the performance output $z(t)$ in (\ref{z}),  should be selected such that the optimization problem  (\ref{lmi}) is feasible. Moreover,  we should prevent rapid   oscillations in the output of the system. If we choose smaller values for the parameters $\gamma_s$ and $\gamma_v$,  the system oscillations increase, while the system response converges more slowly that is not desirable. As a result, there is a trade-off between the convergence rate of the system response and its quality in terms of the amplitude of the oscillations. Accordingly,  adjusting the values of  the parameters $\gamma_s$ and $\gamma_v$ is an important part of the controller design.    %As a result, there exists a trade-off in choosing appropriate values for $\gamma_s$ and $\gamma_v$.
 %Also, smaller values result as slow response of
%the controller. To keep the system response fast, the state
%weights  are chosen to be
In order to have an appropriate output behaviour of the traffic flow system, we choose $\gamma_s=0.03$, $\gamma_v=0.15$, and $\gamma_u=1$.
%The control effort weight is kept unity $\gamma_u=1$.

As the next experiment, we assume that, at $t=20 \:s$, the 7\textsuperscript{th} vehicle is decelerated at $-3\:\mathrm{m}/\mathrm{s}^2$ for $3 \:s$ (this perturbation can be due to road bottlenecks). In Figs. \ref{HDV} (a,b), the trajectory and the velocity profile of all vehicles when there is no CAV in the system are illustrated. One can see  in these figures that  the perturbation does not vanish, and a nonlinear wave appears that propagates against the traffic flow.
On the other hand, when one CAV is added to the traffic system and is controlled by the proposed strategy in Section \ref{se}, one can see in Figs. \ref{cav} (a,b) that  the stop-and-go-wave  can be quickly dissipated, and  the traffic flow is stabilized  to the equilibrium point.

% In this case, from Fig.~\ref{pert}(b), we observe that the perturbations are dampened, and the traffic flow is stabilized to the average velocity. It should be noted that employing the proposed controller in this paper leads to a  higher convergence rate of velocity to the average value as compared to the structured optimal control method in \cite{wang2019controllability}. In \cite{wang2019controllability}, the system is stabilized to the average velocity within 40 s, while here, the convergence occurs within 10 s. %Moreover, in the transient phase, we have less oscillations by using an output-feedback controller.
%We can also change the desired equilibrium velocity to 16 and 14 $\mathrm{m}/\mathrm{s}$, respectively. Considering Fig.~\ref{pert}(c)--(d), one can see that in these cases, the single CAV is still capable of stabilizing traffic flow and steering it towards the new equilibrium points, which verifies the stabilizability of mixed traffic system with only one CAV.

%As the next experiment, we consider a case where, at $t=20 \:s$, the 7th vehicle is decelerated at $-5\:\mathrm{m}/\mathrm{s}^2$ for $5 \:s$. This perturbation does not vanish when there is no CAV in the platoon (Fig. \ref{HDV} (a)), and a nonlinear wave appears that propagates against the traffic flow (Fig. \ref{HDV} (a), (b)).  On the other hand, when one CAV is added to the traffic system and is controlled by the proposed strategy, one can see in Fig. \ref{HDV} (a), (b) that  the stop-and-go-wave  can be quickly dissipated, and  the traffic flow is stabilized  to the equilibrium point.

\begin{figure}[t]

\includegraphics[width=.5\textwidth]{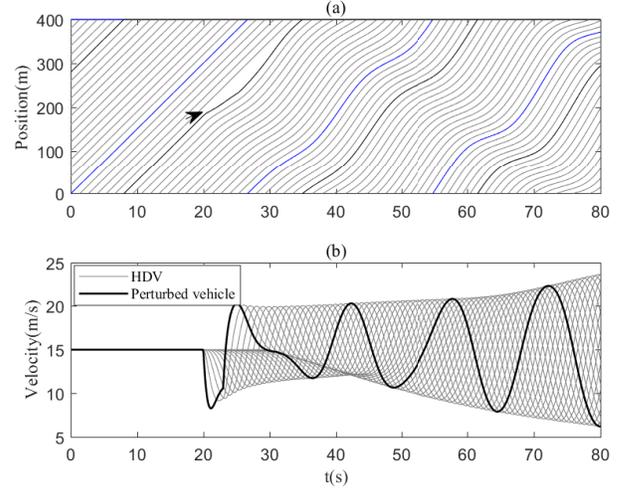}
\centering
\caption{(a) Trajectory of all vehicles (the trajectory of the perturbed vehicle and vehicle no. 20 are, respectively,  shown by black and  blue lines) when there is no CAV in the platoon (a disturbance is added to the acceleration of vehicle 7 at $t=20\:s$);  (b) Velocity profile of vehicles.}
\label{HDV}
\end{figure}

\begin{figure}[t]

\includegraphics[width=.5\textwidth]{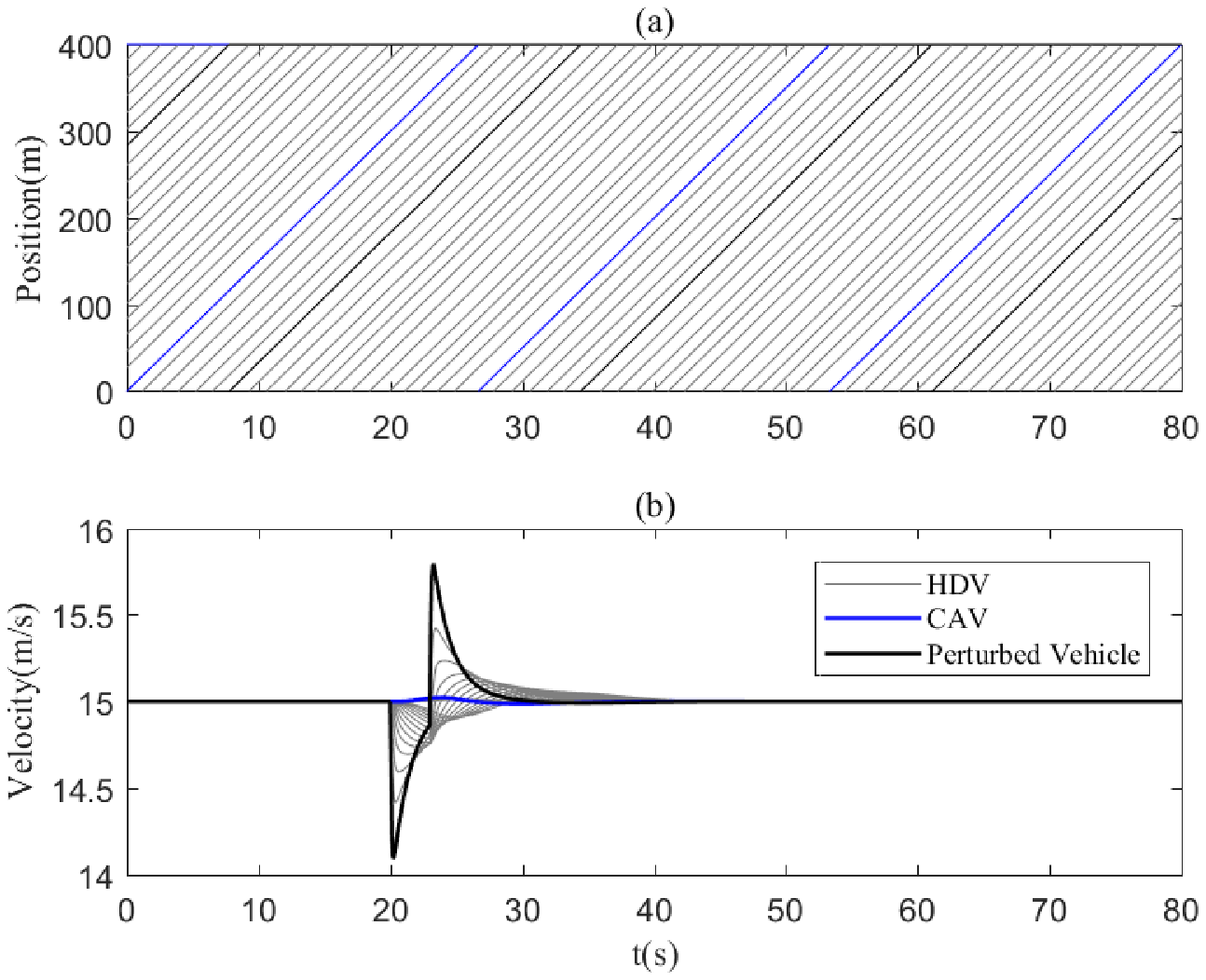}
\centering
\caption{(a) Trajectory of all vehicles (the trajectory of the perturbed vehicle and vehicle no. 20 are, respectively,  shown by black and  blue lines) when there is one CAV under the proposed control is Section \ref{se} (a disturbance is added to the acceleration of vehicle 7 at $t=20\:s$); (b) Velocity profile of vehicles.}
\label{cav}
\end{figure}

\subsection{Comparison with Existing Results}
In this part, we  compare the efficiency of the output-feedback controller proposed in Section \ref{se} with some of  the existing control strategies in the literature. In particular, we compare our result with an optimal control strategy proposed in \cite{wang2020controllability} and  the two heuristic control methods presented in \cite{stern2018dissipation}, that is, Follower-Stopper and PI with Saturation.

At any stage of the new experiment, we assume that one of the HDVs is decelerated at $-3\:\mathrm{m}/\mathrm{s}^2$ for $3 \:s$. Thus, the disturbance signal is a pulse with a time duration of $3 \:s$ and an amplitude of $ -3\:\mathrm{m}/\mathrm{s}^2$.  The initial velocity of all vehicles is also  $15\:\mathrm{m}/\mathrm{s}$.

In Fig.~\ref{cav0}~(a), after applying different control strategies, we illustrate the maximum absolute value of the spacing error during the overall process with respect to the index of a single perturbed HDV. Then, one can observe that with our proposed controller,  the absolute value of the spacing error is much less   than the two heuristic methods \cite{stern2018dissipation} and the optimal control strategy  in \cite{wang2020controllability}. In fact, with our controller, the spacing remains close to the equilibrium spacing over the whole process, and its fluctuations around the desired spacing  are considerably smaller compared to the other strategies. %Moreover, it leaves no large gap from the preceding vehicle

An another performance measure, in Fig.\ref{cav0}~(b), associated with any control strategy,  a linear quadratic cost that is defined as the energy of the performance output $z(t)$ in (\ref{z}) with respect to the different position of the perturbation is illustrated. As is evident in this figure, thanks to considering a system-level performance, our control methodology not only leads to a much better efficiency compared to the heuristic strategies in \cite{stern2018dissipation}, but it has also a substantially smaller cost than the optimal controller proposed  in \cite{wang2020controllability}. In fact,  in \cite{wang2020controllability}, a sub-optimal solution for synthesizing a static controller has been developed, while in this work, a global optimal dynamic controller with more degrees of freedom is designed.

\begin{figure}[t]
%\raggedleft
\includegraphics[width=.5\textwidth]{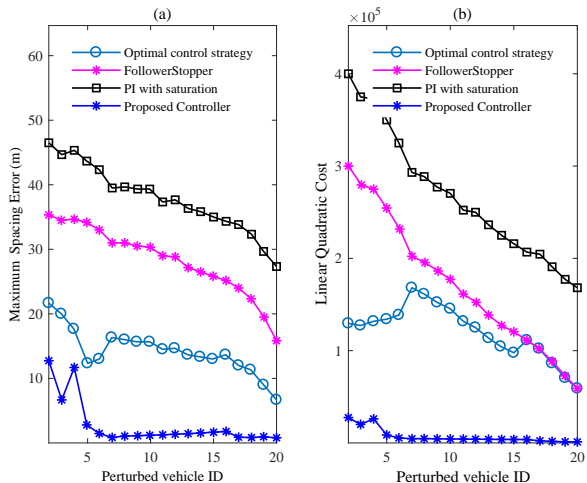}
%\raggedleft
\caption{Comparison of the results for four different control strategies; (a)  $\max_{t} |s_1(t)-s_1^*|$ w.r.t. the index of the perturbed HDV; (b) Energy of the performance output, i.e., $||z(t)||_{2}=\int_{t=0}^{\infty}x^T(t)\mathcal{T}x(t)+u^T(t)\mathcal{Q}u(t)$, w.r.t. the index of the perturbed HDV.   }
\label{cav0}
\end{figure}

 \subsection{Robustness against   Disturbances and Uncertainties}

 As the last experiment, we assume that some parameters appearing in the dynamic model of   HDVs are uncertain. However, in our case study, the nominal values of the parameters  are the same for all HDVs. In fact, instead of considering heterogeneous HDVs, we assume that there are a number of vehicles with homogeneous nominal dynamic models that may also include uncertainties. In order to dampen the perturbations of the mixed traffic system in the presence of parametric uncertainties, we design  an output-feedback controller, using the procedure proposed in Section \ref{kon2}.

 In this experiment, we assume that the parameters $\alpha_i$, $\theta_i$, and $s_{i,go}$ (see Section \ref{hii}) are unknown.  The uncertain parameters are
distributed  around the nominal values $\alpha_{\mathcal{N}}=0.6$, $\theta_{\mathcal{N}}=0.9$, and $s_{\mathcal{N},go}=35$. Then, one can define $\alpha_i=\alpha_{\mathcal{N}}+\Delta \alpha$, $\theta_i=\theta_{\mathcal{N}}+\Delta \theta$, $s_{i,go}=s_{\mathcal{N},go}+\Delta s_{go}$, where $-0.1\leq \Delta \alpha\leq 0.1$, $-0.1\leq \Delta \theta\leq 0.1$, and $-5\leq \Delta s_{go}\leq 5$.

%The following equations represent the random variation of the parameters
%\begin{equation*}
%\begin{aligned}
%&\alpha_i=\alpha_{ni}+\frac{\alpha_{i,max}-\alpha_{i,min}}{2}rand(t)\\ &\theta_i=\theta_{ni}+\frac{\theta_{i,max}-\theta_{i,min}}{2}rand(t)\\ &s_{i,go}=s_{ni,go}+\frac{s_{imax,go}-s_{imin,go}}{2}rand(t)\end{aligned}
%\end{equation*}

%where $\alpha_{i,max}=0.7,\alpha_{i,min}=0.5, \theta_{i,max}=1, \theta_{i,min}=0.8, s_{imax,go}=40, s_{imin,go}=30$.
%rand(t) is a stochastic process with uniform distribution in the
%interval $[-1,1]$ for any given time t.

We first assume in this experiment that the  initial velocity of the vehicles is randomly chosen from the distribution $15+U[-4,4] \:\mathrm{m}/\mathrm{s}$.
We compute the controller parameters associated with  three different   equilibrium velocities, that is,  $14$, $15$, and $16\:\mathrm{m}/\mathrm{s}$. As observed in Figs.~\ref{pert}~(a--c), in these three  cases,  by using the control strategy proposed in Section \ref{kon2}, the perturbations occurring in the traffic flow system are dampened, and the velocities of all vehicles converge to the equilibrium points within a short time.

Next, we assume that, at $t=20 \:s$, the  vehicle no. 7 brakes at $-3\:\mathrm{m}/\mathrm{s}^2$ for $3 \:s$, and the single CAV is under the proposed control in Section \ref{kon2}. In Fig.~\ref{pertt1}~(a,b), the velocity profile of all vehicles from a 3-dimensional (3D) and 2-dimensional (2D) perspective are illustrated. One can see the with the proposed output-feedback controller, the perturbation is attenuated very quickly in this traffic system. Moreover, as observed in Fig.~\ref{pertt2}, the spacing of HDVs and the CAV from the preceding vehicles converges to the equilibrium value in a few seconds.

\begin{figure}[t]
\includegraphics[width=.5\textwidth]{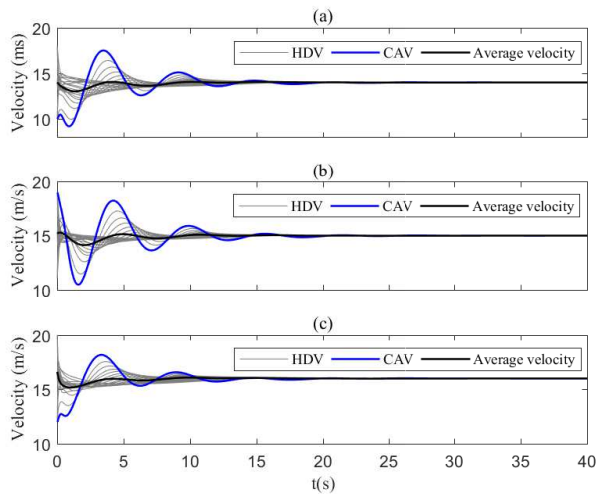}
\centering
\caption{ (a)--(c) Velocity profile of vehicles when there is one CAV under the proposed control in Section \ref{kon2} in the presence of disturbances and the parametric uncertainties, and $v^*$ is 14, 15, and 16 $\mathrm{m}/\mathrm{s}$ for (a), (b), and (c), respectively.}
\label{pert}
\end{figure}

\begin{figure}[t]
\includegraphics[width=.5\textwidth]{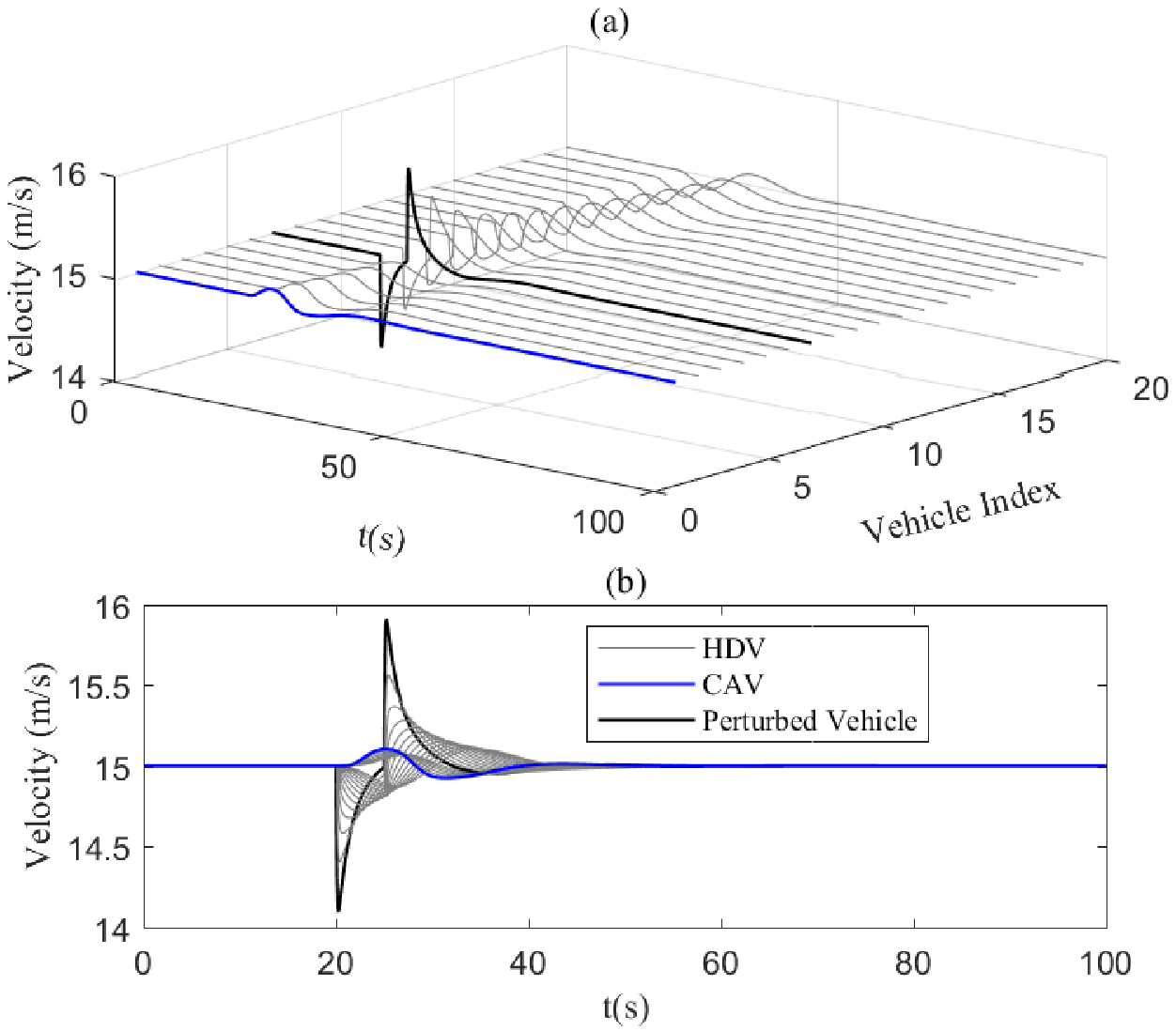}
\centering
\caption{ (a)--(b) Velocity profile of all vehicles over time from 3D and 2D perspectives in the presence of disturbances and parametric uncertainties when the single CAV is under the proposed control in Section \ref{kon2} (the blue (resp., black) line represents the spacing of the CAV (resp., the perturbed vehicle), i.e. $s_1(t)$ (resp., $s_7(t)$), and the grey lines denote the spacing of other HDVs).  }
\label{pertt1}
\end{figure}

\begin{figure}[t]
\includegraphics[width=.5\textwidth]{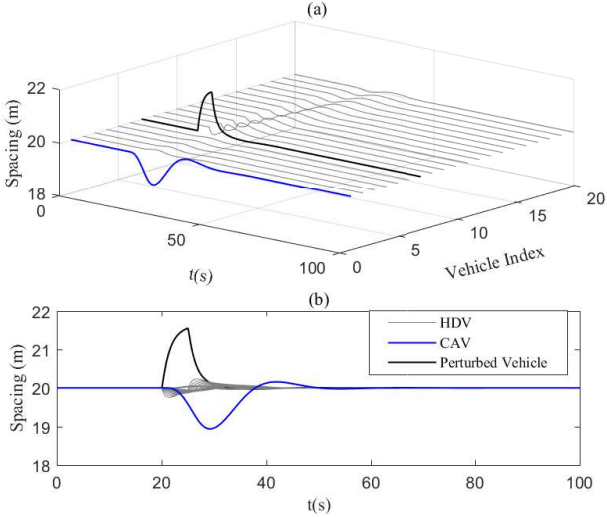}
\centering
\caption{ (a)--(b) Spacing profile of all vehicles over time from 3D and 2D perspective in the presence of disturbances and parametric uncertainties when the single CAV is under the proposed control in Section \ref{kon2} (the blue (resp., black) line represents the spacing of the CAV (resp., the perturbed vehicle), i.e. $s_1(t)$ (resp., $s_7(t)$), and the grey lines denote the spacing of other HDVs). }
\label{pertt2}
\end{figure}

\section{Conclusion}

In this work, the stabilizability and detectability of a mixed traffic  system along a ring-road and an open-road has been studied and analyzed. It has been shown that in both cases, the system under investigation is stabilizable when there is one single CAV in the platoon. Moreover, we have demonstrated that the system is detectable even if the state information of only one vehicle is directly measured. Furthermore, by considering limited communication ability of the CAV to receive state information from its neighboring vehicles, an $H_{\infty}$ dynamic output-feedback controller has been designed, that can provably smoothen the traffic flow even in the presences of abrupt and large disturbances (e.g.\ sharp deceleration downstream). We have also presented an $H_{\infty}$ control method in the presence of parametric uncertainties in the dynamics of HDVs.  The effectiveness of the proposed methods to achieve traffic flow stability has been verified and demonstrated through numerical simulation experiments. Finally, the numerical results are compared to another approach from the literature. Future work is going to deal with the generalization of these results to more complex traffic patterns.

\appendix
We first present two results that will
be used in the proof of Theorem \ref{uncertain} and then discuss the proof.

\begin{lem}{\cite{wang1992robust}}
(Young inequality) Given
 matrices $D, F$, and $S$ of  appropriate dimensions, the inequality
$DFS +(DFS)^{T}\preceq \epsilon^{-1}DD^T+\epsilon S^TS$ holds for some  scalar $\epsilon>0$  if we have $F^T F\preceq I$.\label{young}
\end{lem}

\begin{lem}{\cite{zhang2006schur}}
 (Schur complement) Consider the
matrices  $W_1$, $W_2$, and $W_3$ with the appropriate dimensions,  where
$W_1=W_1^T\succ0$. Then,  the matrix inequality $W_1+W_3^T W_2^{-1} W_3\prec 0$ is equivalent to
$$\begin{bmatrix}
W_1 & W_3^T\\ W_3 & -W_2 \end{bmatrix}\prec 0.$$\label{schur}
\end{lem}

\emph{Proof of Theorem \ref{uncertain}:}
As shown in section (\ref{se}), based on the bounded real lemma (BRL), if the matrix inequality (\ref{si3}) is satisfied, then $||T_{zd}||_{\infty}<\gamma$.
Incorporating $A$ from (\ref{delta}) to (\ref{si3}), we get
$$\Omega_4=\Omega_3+\Delta \Omega\prec0,$$ where \begin{equation*}
\begin{aligned}
\Omega_3=\begin{bmatrix}
\Omega_{11}&\Omega_{12}&\Omega_{13}&\Omega_{14}\\\Omega_{21}&\Omega_{22}&\Omega_{23}&\Omega_{24}\\ \Omega_{31}&\Omega_{32}& \Omega_{33}&\Omega_{34}
\\\Omega_{41}&\Omega_{42}& \Omega_{43}&\Omega_{44}\end{bmatrix} \label{sii4}\end{aligned},\end{equation*}
with
\begin{equation*}
\begin{aligned}
        &\Omega_{11}=
A_{\mathcal{N}}X+XA_{\mathcal{N}}^T+BC_k M^T+ MC_k^T B^T,\\&\Omega_{12}=\Omega^T_{21}=MA_k^TN^T+XC^TB_k^TN^T+MC_k^TB^TY,\\&\quad\quad\quad\quad+XA_{\mathcal{N}}^TY+A_{\mathcal{N}}\\&\Omega_{22}=A_{\mathcal{N}}^T Y+YA_{\mathcal{N}}+NB_kC+C^T B_k^T N^T,\\&\Omega_{13}=\Omega^T_{31}=B_d, \;\; \Omega_{14}=\Omega^T_{41} =XC_z^T+MC_k^TD^T_z, \\&\Omega_{23}=\Omega^T_{32}=YB_d, \;\;\Omega_{24}=\Omega^T_{42}=C_z^T, \;\;\Omega_{33}=-\gamma^2I,\\& \Omega_{34}=\Omega_{43}=0,\;\;\Omega_{44}=- I,\\&\\\mbox{and}\\&\Delta\Omega=\begin{bmatrix}
\Delta AX+X\Delta A^T & \Delta A+X\Delta A^T Y & 0 & 0\\\Delta A^T+Y\Delta A X&\Delta A^T Y+Y\Delta A &0&0\\ 0&0& 0&0\\0&0&0&0
\end{bmatrix}.\end{aligned}
\end{equation*}
By substituting $\Delta A=LFR$, and  applying
condition (\ref{conun}) and Lemma \ref{young}, $ \Omega_4$  can be bounded as
\begin{equation*} \Omega_4 \preceq \Omega_3+\Sigma_{i=1}^4 \Upsilon_i\prec0,\end{equation*}
where
\begin{equation*}
\begin{aligned}
&\Upsilon_{1}=\begin{bmatrix}
L\\0\\0\\0\end{bmatrix} \epsilon_1 \begin{bmatrix}
L^T&0&0&0\end{bmatrix}+\begin{bmatrix}
XR^T\\0\\0\\0\end{bmatrix} \epsilon_1^{-1} \begin{bmatrix}
RX&0&0&0\end{bmatrix},\\&\Upsilon_{2}=\begin{bmatrix}
L\\0\\0\\0\end{bmatrix} \epsilon_2 \begin{bmatrix}
L^T&0&0&0\end{bmatrix}+\begin{bmatrix}
0\\R^T\\0\\0\end{bmatrix} \epsilon_2^{-1} \begin{bmatrix}
0&R&0&0\end{bmatrix},\\&\Upsilon_{3}=\begin{bmatrix}
XR^T\\0\\0\\0\end{bmatrix} \delta \begin{bmatrix}
RX&0&0&0\end{bmatrix}+\begin{bmatrix}
0\\YL\\0\\0\end{bmatrix} \delta^{-1}\begin{bmatrix}
0&L^TY&0&0\end{bmatrix}, \\&\Upsilon_{4}=\begin{bmatrix}
0\\R^T\\0\\0\end{bmatrix} \epsilon_3 \begin{bmatrix}
0&R&0&0\end{bmatrix}+\begin{bmatrix}
0\\YL\\0\\0\end{bmatrix} \epsilon_3^{-1} \begin{bmatrix}
0&L^TY&0&0\end{bmatrix}.  \label{si44}\end{aligned} \end{equation*}

Now, by defining $
\hat{A}=NA_kM^T+NB_k CX+YBC_k M^T+YA_{\mathcal{N}}X$, $ \hat{B}=NB_k $, $\hat{C}=C_k M^T$, $\eta=\gamma^2$ and selecting $\delta=1$, we obtain
$$\Omega_5+\Gamma_{12}\Gamma^{-1}_{22}\Gamma^T_{12}\prec0$$
where
\begin{equation*}
\begin{aligned}
\Omega_5=\begin{bmatrix}
\Omega_{11}&\Omega_{12}&\Omega_{13}&\Omega_{14}\\\Omega_{21}&\Omega_{22}&\Omega_{23}&\Omega_{24}\\ \Omega_{31}&\Omega_{32}& \Omega_{33}&\Omega_{34}
\\\Omega_{41}&\Omega_{42}& \Omega_{43}&\Omega_{44}\end{bmatrix}, \label{siii4}\end{aligned}\end{equation*}
with
\begin{equation*}
\begin{aligned}
        &\Omega_{11}=
A_{\mathcal{N}}X+XA_{\mathcal{N}}^T+B\hat{C}+ \hat{C}^T B^T+(\epsilon_1+\epsilon_2)LL^T,\\&\Omega_{12}=\Omega^T_{21}=\hat{A}^T+A_{\mathcal{N}}\\&\Omega_{22}=A_{\mathcal{N}}^T Y+YA_{\mathcal{N}}+\hat{B}C+C^T \hat{B}^T+\epsilon_3R^TR,\\&\Omega_{13}=\Omega^T_{31}=B_d, \;\; \Omega_{14}=\Omega^T_{41} =XC_z^T+\hat{C}^TD^T_z, \\&\Omega_{23}=\Omega^T_{32}=YB_d, \;\;\Omega_{24}=\Omega^T_{42}=C_z^T, \;\;\Omega_{33}=-\eta I,\\& \Omega_{34}=\Omega_{43}=0,\;\;\Omega_{44}=- I
.\end{aligned}
\end{equation*}
%
%\vspace{-5mm}
Moreover, we have
\begin{equation*}
\begin{aligned}
\Gamma_{12}=&\begin{bmatrix}
XR^T&0&0&YL&0\\0&R^T&YL&0&XR^T\\0&0&0&0&0\\0&0&0&0&0
\end{bmatrix},\\
\Gamma_{22}=&-\mathrm{diag}(\epsilon_1I,\epsilon_2I,\epsilon_3I,I,I).
\end{aligned}
\end{equation*}
Finally, by applying Lemma \ref{schur}, the LMIs in  (\ref{lmiun}) are obtained.

\section*{Acknowledgement}

This work was supported in part by the Swiss National Science Foundation (SNSF) under the project RECCE, ``Real-time traffic estimation and control in a connected environment'',
contract No. 200021-188622.

\bibliographystyle{IEEEtran}
\bibliography{library}

\end{document}